\documentclass{aims}

\usepackage{amsmath}
\usepackage{paralist}
\usepackage{graphics}
\usepackage{epsfig}
\usepackage[colorlinks=true]{hyperref}
\hypersetup{urlcolor=blue, citecolor=red}

\def\currentvolume{X}
\def\currentissue{X}
\def\currentyear{200X}
\def\currentmonth{XX}
\def\ppages{X--XX}

\newtheorem{theorem}{Theorem}[section]

\theoremstyle{definition}

\allowdisplaybreaks[4]

\newcommand{\C}[1]{\mathbf{C}^{#1}}

\newcommand{\norma}[1]{{\left\|#1\right\|}}
\newcommand{\reali}{{\mathbb{R}}}

\newcommand{\BV}{\mathbf{BV}}

\renewcommand{\epsilon}{\varepsilon}
\renewcommand{\phi}{\varphi}
\renewcommand{\L}[1]{{\mathbf{L}^#1}}
\newcommand{\W}[2]{{\mathbf{W}^{#1,#2}}}
\renewcommand{\O}{\mathcal{O}(1)}
\newcommand{\Lloc}[1]{{\mathbf{L}_{loc}^{#1}}}
\newcommand{\tv}{\mathinner{\rm TV}}
\newcommand{\caratt}[1]{{\chi_{\strut#1}}}

\renewcommand{\div}{\mathinner{\rm div}}

\newcommand{\pt}{\partial}

\newcommand{\Lip}{\mathinner\mathbf{Lip}}
\newcommand{\rpic}{{\mathbb{R}^+}}
\newcommand{\Cloc}[1]{{\mathbf{C}_{loc}^{#1}}}

\renewcommand{\d}[1]{\mathinner{\mathrm{d}{#1}}}

\title[Conservation Laws in Crowds Modeling]{Conservation Laws\\ in the
  Modeling of Moving Crowds}

\author[R.M.Colombo, M.Garavello, M.L\'ecureux-Mercier, N.Pogodaev]{}


\email{rinaldo@ing.unibs.it}
\email{mauro.garavello@unimib.it}
\email{magali.lecureux-mercier@univ-orleans.fr}
\email{nickpogo@gmail.com}

\begin{document}
\maketitle

\centerline{Rinaldo M.~Colombo} \medskip
{\footnotesize\centerline{Department of Mathematics, Brescia
    University, Brescia, Italy}}

\bigskip

\centerline{Mauro Garavello} \medskip
{\footnotesize\centerline{Department of Mathematics and Applications,
    Milano -- Bicocca University, Milano, Italy}}

\bigskip

\centerline{Magali L\'ecureux-Mercier} \medskip
{\footnotesize\centerline{Technion, Israel Institute of Technology,
    Haifa, Israel}}

\bigskip

\centerline{Nikolay Pogodaev} \medskip
{\footnotesize\centerline{Russian Academy of Sciences, Irkutsk,
    Russia}}

\begin{abstract}
  Models for crowd dynamics are presented and compared. Well posedness
  results allow to exhibit the existence of optimal controls in
  various situations. A new approach not based on partial differential
  equations is also briefly considered.
\end{abstract}

\section{Introduction}
\label{sec:Intro}

\par From a macroscopic viewpoint, a moving crowd can be described
through its density $\rho = \rho (t,x)$, a function of time $t \in
\reali^+$ and space $x \in \reali^2$ attaing values in $[0,1]$. In
standard situations, the number of pedestrians is conserved, so that
$\int_{\reali^2} \rho (t,x) \d{x}$ is independent of $t$. Hence, it is
natural to use the \emph{conservation law}
\begin{equation}
  \label{eq:1}
  \partial_t \rho + \div_x (\rho \,  V) = 0 \,.
\end{equation}
Any model of this kind depends on the \emph{speed law} that defines
the velocity $V$ of the crowd as a function of $t$, $x$, $\rho$,
$\ldots$ A simple version of~\eqref{eq:1} is obtained assigning
\begin{equation}
  \label{eq:2}
  V = v (\rho) \, \vec{v}(x)
  \quad \mbox{ with } \quad
  \begin{array}{l}
    v \in \C2 ([0,1]; \reali^+)
    \mbox{ non increasing and } v (1)=0 \,,
    \\
    \vec{v} \in \C2 (\reali^2; \mathbb{S}^1) \,.
  \end{array}
\end{equation}
In this case, Kru\v zkov Theorem~\cite[Theorem~1]{Kruzkov} applies and
ensures that the Cauchy problem for~\eqref{eq:1}--\eqref{eq:2} has a
unique solution in $\C{0,1} \left(\reali^+; \L1 (\reali^2;
  [0,1])\right)$ which depends Lipschitz continuously from the data
and, by~\cite[Theorem~2.6]{ColomboMercierRosini}, also from $v$ and
$\vec{v}$.

According to~\eqref{eq:2}, at time $t$ the pedestrian at $x$ moves
along a prescribed trajectory, an integral curve of $\vec{v}$, with a
speed $v (\rho)$ that depends on $\rho$ evaluated at point $x$ and
time $t$. On the contrary, Section~\ref{sec:S} is devoted
to~\eqref{eq:1} with the speed of the individual at $x$ depending on
an \emph{average} of the density $\rho$ in a \emph{neighborhood} of
$x$. The resulting model has a rich analytical structure, the
solutions being also \emph{differentiable} with respect to the data
and to the speed law.

In Section~\ref{sec:R} the direction chosen by the pedestrian at $x$
depends from an \emph{average} of the density \emph{gradient} $\nabla
\rho$ around $x$, while his/her speed depends from $\rho$ evaluated at
$x$. The resulting solutions display qualitative properties usually
seen in context where individuals have a proper volume such as the
\emph{Braess paradox}~\cite{BraessParadox} and the formation of
queues~\cite{HelbingEtAlii2001}.

If the various individuals have different destinations then it is
possible to subdivide the crowd under consideration into different,
say $n$, populations with densities $\rho_1, \ldots, \rho_n$, each
having a different destination. The resulting model
\begin{equation}
  \label{eq:3}
  \partial_t \rho_i + \div_x (\rho_i \, V_i) = 0
  \qquad i=1, \ldots, n
\end{equation}
consists of a system of conservation laws that, when $n=1$, reduces
to~\eqref{eq:1}. The results in both Section~\ref{sec:S} and
Section~\ref{sec:R} can be extended to this more general setting.

Finally, Section~\ref{sec:D} approaches the problem of driving a crowd
with a few moving individuals. First, a model based on~\eqref{eq:1} is
recalled and then an approach based on differential inclusions is
presented. The latter approach, developed following~\cite{BressanFire,
  BressanZhang}, neglects the crowd internal dynamics and allows for a
simpler analytical framework.

We refer for instance to~\cite{BellomoDogbe2011} for an account of the
fast development of the recent macroscopic modeling of crowd
dynamics. Moreover, measure valued conservation laws were considered
in~\cite{CristianiPiccoliTosin, PiccoliTosin}; the results
in~\cite{MauryEtAl} deal with constrained velocity models; various 1D
attempts are found in~\cite{AmadoriDiFrancesco, ColomboRosini2005,
  ColomboRosini2009, DiFrancescoMarkowichWolfram,
  GoatinRosini}. Throughout, for the basic results in the theory of
conservation laws we refer to~\cite{BressanLectureNotes,
  DafermosBook}.

\section{NonLocal Speed Choice}
\label{sec:S}

Consider~\eqref{eq:1} with the nonlocal speed law
\begin{equation}
  \label{eq:4}
  V (\rho) = v\left(\rho (t) * \eta\right) \, \vec{v} \,.
\end{equation}
Here, the speed $v$ at time $t$ of the pedestrian at $x$ depends on
the averaged density $\left(\rho(t)*\eta\right)(x) = \int_{\reali^2}
\rho(t, x-y) \, \eta(y) \, \d{y}$. The direction of the velocity is
given by the (fixed) vector $\vec{v} (x)$.

For simplicity, we state the results below in $\reali^2$. However, the
case where the region available to the crowd is constrained by, say,
walls or doors can be easily recovered in the present framework, along
the technique used in~\cite{ColomboGaravelloLecureux,
  ColomboGaravelloLecureuxM3AS}

As is typical whenever Kru\v zkov techniques apply, space dimension
$2$ plays no role and the results below can be extend to $\reali^n$.

Existence and uniqueness of a solution to the Cauchy problem
for~\eqref{eq:1}--\eqref{eq:4} follow from the next result.

\begin{theorem}{\cite[Proposition~4.1]{ColomboHertyMercier},
    \cite[Theorem~2.2]{ColomboLecureuxPerDafermos}}
  \label{thm:panicK}
  Let $v \in(\C2\cap \W2 \infty)(\reali; \reali)$, $\vec \nu \in
  (\C2\cap\W21)(\reali^2; \reali^2)$, $\eta \in (\C2\cap\W2
  \infty)(\reali^2; \reali)$. Assume $\rho_o \in (\L1\cap\L
  \infty\cap\BV)(\reali^2; \rpic)$. Then, \eqref{eq:1}--\eqref{eq:4}
  with initial condition $\rho_o$ admits a unique weak entropy
  solution $\rho \in \C0 \left(\rpic; \L1(\reali^2; \rpic) \right)$.
  Furthermore, we have the estimate $\displaystyle
  \norma{\rho(t)}_{\L\infty}\leq \norma{\rho_o}_{\L\infty} e^{Ct}$,
  where the constant $C$ depends on $v$, $\vec \nu$ and $\eta$.
\end{theorem}

The definition of weak entropy solutions is based on Kru\v zkov
notion~\cite[Definition~1]{Kruzkov}, see
also~\cite{ColomboHertyMercier, ColomboLecureuxPerDafermos}. The proof
relies on a contraction argument based on the key estimates provided
by~\cite[Theorem~2.6]{ColomboMercierRosini}.

Another contraction argument, based on tools from optimal transport
theory, allows to extend the above result to the measure valued
setting in~\cite{CrippaMercier}. (Below, $\mathcal{M}^+ (\reali^+)$ is
the set of positive Radon measures on $\reali^2$).

\begin{theorem}{\cite[Theorem~1.1]{CrippaMercier}}
  \label{thm:panicOT}
  Assume $v \in (\L \infty\cap\Lip)(\reali; \reali)$, $\vec \nu \in
  (\L \infty\cap\Lip)(\reali^2; \reali^2)$, $\eta \in (\L
  \infty\cap\Lip)(\reali^2; \rpic)$. Let $\rho_o \in \mathcal{M}^+
  (\reali^2)$. Then, there exists a unique weak measure valued
  solution $\rho \in \L \infty(\rpic; \mathcal{M}^+(\reali^2))$
  to~\eqref{eq:1}--\eqref{eq:4} with initial condition $\rho_o$. If
  furthermore $\rho_o \in \L1(\reali^2; \rpic)$, then $\rho \in
  \C0\left(\reali^+; \L1(\reali^2; \rpic)\right)$.
\end{theorem}

In general, in~\eqref{eq:1}--\eqref{eq:4} no \emph{a priori} uniform
$\L\infty$ bound on the density is possible. Indeed, assume that the
density is $1$ all along the trajectory of the pedestrian at $x$. The
averaged density around $x$ may well be less than $1$, forcing the
pedestrian to proceed and, hence, leading to a increase in the
density. This behavior can be related to the rise of panic,
see~\cite{ColomboRosini2005, ColomboRosini2009}. In the literature,
values of $\rho$ of up to 10 individuals per square meter were
measured, see for instance~\cite{HelbingJohanssonZein}.

Aiming at preventing the insurgence of these phenomena, it is natural
to consider control problems where functionals of the density of the
type
\begin{equation}
  \label{eq:5}
  J_T(\rho_o)
=
\int_0^T \int_{\Omega} f\left(\rho (t,x)\right) \, \d{x} \, \d{t}
  \quad
  \mbox{where }
  \rho \mbox{ solves \eqref{eq:1}--\eqref{eq:4}}
  \mbox{ with datum } \rho_o
\end{equation}
have to be minimized. Here, $\Omega$ is the region where the density
needs to be controlled and $f$ is a $\C1$ function weighing $0$ on
acceptable densities and quickly increasing when $\rho$ approaches
dangerous values. Necessary conditions for the minima of~\eqref{eq:5}
are available once the differentiability of the solution
to~\eqref{eq:1}--\eqref{eq:4} with respect to the initial datum is
proved. This motivates the following result.

\begin{theorem}{\cite[Theorem~4.2]{ColomboHertyMercier}
    \cite[Theorem~2.2]{ColomboLecureuxPerDafermos}}
  \label{thm:Gdiff}
  Let $\rho_o \in (\W2 \infty\cap \W21)(\reali^2; \rpic)$, $r_o \in
  (\W11\cap \L\infty)(\reali^2; \reali)$. Assume $v \in (\C4\cap\W2
  \infty)(\reali;\reali )$, $\vec \nu \in (\C3\cap\W21)(\reali^2;
  \reali^2)$, $\eta \in (\C3\cap \W2 \infty)(\reali^2; \rpic)$.  Then,
  there exists a unique weak entropy solution $r \in \C0 (\reali^+;
  \L1 \left(\reali^2; \reali)\right)$ to the Cauchy problem
  \begin{equation}\label{eq:linearised}
    \pt_t r + \div (r \, v(\rho*\eta) \,\vec \nu(x))
    =
    -\div ( \rho\, v'(\rho*\eta)\,\vec \nu(x))\,,\qquad r(0)=r_o\,.
  \end{equation}
  Furthermore, for all $\rho_o \in (\W21\cap\W2 \infty)(\reali^2;
  \rpic)$ and $r_o \in(\W11\cap\L \infty)(\reali^2; \reali) $, call
  $\rho_h$ the solution to~\eqref{eq:1}--\eqref{eq:4} with initial
  datum $\rho_o + h r_o$. Then, for all $t \in \reali^+$,
  \begin{equation}
    \label{eq:limit}
    \lim_{h\to 0}
    \norma{\frac{\rho_h (t) - \rho (t)}{h} - r (t)}_{\L1}
    = 0
  \end{equation}
  i.e., the solution $\rho$ to~\eqref{eq:1}--\eqref{eq:4} is G\^ateaux
  differentiable in $\rho_o$ along any direction $r_o$.
\end{theorem}

To prove this theorem, first the well posedness
of~\eqref{eq:linearised} is obtained and then the
limit~\eqref{eq:limit} is computed. In both steps, the estimates
in~\cite{ColomboMercierRosini} play a key role. At present, no analog
to Theorem~\ref{thm:Gdiff} is available in the setting of
Theorem~\ref{thm:panicOT}. Indeed, a good definition of G\^ateaux
differentiability on the set of probability measures equipped with the
Wasserstein distance of order 1 is, to our knowledge, not available.

\section{NonLocal Route Choice}
\label{sec:R}

Consider~\eqref{eq:1} with the nonlocal speed law
\begin{equation}
  \label{eq:6}
  V (\rho)
  =
  v (\rho) \, \left(\vec{\nu} (x) + \mathcal{I} (\rho)\right) \,.
\end{equation}
Here, the individual in $x$ at time $t$ moves at the speed
$v\left(\rho (t,x)\right)$ that depends on the density $\rho (t,x)$
evaluated at the same time $t$ and $x$. The vector $\vec{\nu} (x) \in
\reali^2$ is the preferred direction of the pedestrian at $x$, while
$\mathcal{I}(\rho) (x)$ describes how the pedestrian at $x$ deviates
from the preferred direction, given that the crowd distribution is
$\rho$. Thus, the individual at time $t$ in $x$ is assumed to move in
the direction of the vector $\vec{\nu}(x) + \left( \mathcal{I} \left(
    \rho(t) \right) \right) (x)$. The basic well posedness result
for~\eqref{eq:1}--\eqref{eq:6} is the following.

\begin{theorem}{~\cite[Theorem~2.1,
    Theorem~2.2]{ColomboGaravelloLecureuxM3AS}}
  \label{thm:existence}
  Let the following conditions hold:
  \begin{description}
  \item[(v)] $v \in \C2 (\reali; \reali)$ is non increasing, $v (0) =
    V$ and $v (R)=0$ for fixed $V,R>0$.

  \item[($\boldsymbol{\vec{\nu}}$)] $\vec{\nu} \in (\C2 \cap
    \W1\infty) (\reali^2; \reali^2)$ is such that $\div \vec{\nu} \in
    (\W11 \cap \W1\infty) (\reali^2; \reali)$.

  \item[(I)] $\mathcal{I} \in \C0 \left( \L1(\reali^2; [0,R]); \C2
      (\reali^2; \reali^2) \right)$ satisfies the estimates:
    \begin{enumerate}[\bf({I}.1)]
    \item There exists an increasing $C_I \in \Lloc\infty(\rpic;
      \rpic)$ such that, for all $r \in \L1(\reali^2;[0,R])$,
      $\norma{\mathcal{I}(r)}_{\W1\infty} \! \leq
      C_I(\norma{r}_{\L1})$ and $\norma{\mathcal{\div I}(r)}_{\L1}
      \leq C_I(\norma{r}_{\L1})$.

    \item There exists an increasing $C_I \in \Lloc\infty(\rpic;
      \rpic)$ such that, for all $r \in \L1(\reali^2;[0,R])$,
      $\norma{\nabla \div \mathcal{I}(r)}_{\L1} \leq
      C_I(\norma{r}_{\L1})$.

    \item There exists a constant $K_I$ such that for all $r_1, r_2\in
      \L1(\reali^2; [0,R])$,
      \begin{eqnarray*}
        \norma{\mathcal{I}(r_1) - \mathcal{I}(r_2)}_{\L\infty}
        & \leq &
        K_I \cdot \norma{r_1-r_2}_{\L1} \,,
        \\
        \norma{\mathcal{I}(r_1) - \mathcal{I}(r_2)}_{\L1}
        +
        \norma{\div(\mathcal{I}(r_1 ) - \mathcal{I}(r_2))}_{\L1}
        &\leq &
        K_I \cdot \norma{r_1-r_2}_{\L1} \,.
      \end{eqnarray*}
    \end{enumerate}
  \end{description}
  \noindent Choose any $\rho_o \in (\L1 \cap \BV)
  (\reali^2;[0,R])$. Then, there exists a unique weak entropy solution
  $\rho \in \C0 \left(\rpic; \L1(\reali^2; [0,R]) \right)$
  to~\eqref{eq:1}--\eqref{eq:6}. Moreover, $\rho$ satisfies the bounds
  \begin{eqnarray*}
    \norma{\rho(t)}_{\L1}
    & = &
    \norma{\rho_o}_{\L1}\,,  \mbox{ for a.e. } t\in \rpic\,,
    \\
    \tv(\rho(t))
    & \leq &
    \tv(\rho_o) \, e^{k t}
    + \frac{\pi}{4} te^{k t} N \norma{q}_{\L\infty([0,R])}
    \left(
      \norma{\nabla\div \vec{\nu}}_{\L1}
      +
      C_I(\norma{\rho_o}_{\L1})
    \right)\,,
  \end{eqnarray*}
  where $k= (2N+1) \norma{q'}_{\L\infty([0,R])} \left( \norma{\nabla
      \vec{\nu}}_{\L\infty} + C_I(\norma{\rho_o}_{\L1}) \right)$.  If
  also the speed law
  \begin{equation}
    \label{eq:2eq}
    V' (\rho)
    =
    v' (\rho) \, \left(\vec{\nu'} (x) + \mathcal{I}' (\rho)\right)
  \end{equation}
  satisfies the same assumptions, then the solution $\rho$
  to~\eqref{eq:1}--\eqref{eq:6} and $\rho'$
  to~\eqref{eq:1}--\eqref{eq:2eq}, with data $\rho_o,\rho_o' \in (\L1
  \cap \BV) (\reali^2;[0,R])$, satisfy
  \begin{eqnarray*}
    \norma{\rho_1 (t) - \rho_2 (t)}_{\L1}
    \! & \leq & \!
    \left(1 + C (t)\right) \!
    \norma{\rho_{0,1}-\rho_{0,2}}_{\L1} \!
    % \\
    % & + &
    +
    C (t) \!
    \left(
      \norma{q_1 - q_2}_{\W1\infty}\!
      +
      d\left(\mathcal{I}_1,\mathcal{I}_2\right)
    \right)
    \\
    & & +
    C (t)
    \left(
      \norma{\vec{\nu}_1 - \vec{\nu}_2}_{\L\infty}  +
      \norma{\div(\vec{\nu}_1 - \vec{\nu}_2)}_{\L1}
    \right)
  \end{eqnarray*}
  where
  \begin{displaymath}
    d (\mathcal{I}_1, \mathcal{I}_2)
    \!=
    \sup
    \left\{
      \norma{\mathcal{I}_1 (\rho) -\mathcal{I}_2(\rho)}_{\L\infty} \!
      +
      \norma{\div \left(\mathcal{I}_1 (\rho) - \mathcal{I}_2(\rho)\right)}_{\L1}
      \! \colon \! \rho \in \L1 (\reali^2; [0,R])
    \right\} .
  \end{displaymath}
  The map $C \in \C0 (\reali^+; \reali^+)$ vanishes at $t = 0$ and
  depends on $\tv(\rho_{0,1})$, $\norma{\rho_{0,1}}_{\L1}$,
  $\norma{\vec{\nu}_1}_{\L\infty}$, $\norma{\div \vec{\nu}_1}_{\W11}$,
  $\norma{q_1}_{\W1\infty}$, $\norma{q_2}_{\W1\infty}$.
\end{theorem}

% \subsection{Braess paradox}
% \label{sse:braess}

In operation research, Braess paradox states that adding extra
capacity to a network can, in some cases, reduce the overall
performance of the network, see~\cite{BraessParadox}. A relevant
problem in the design of escape routes is the planning of suitable
devices that reduce the exit time. The
model~\eqref{eq:1}--\eqref{eq:6} allows to show that the careful
introduction of suitable obstacles in suitable locations does indeed
reduce the exit time.  In fact, these obstacles reduce congested areas
at the sides of the door jambs.

\begin{figure}[htpb]
  \begin{minipage}[c]{0.43\linewidth}
    \includegraphics[width=\textwidth, trim=60 120 50
    120]{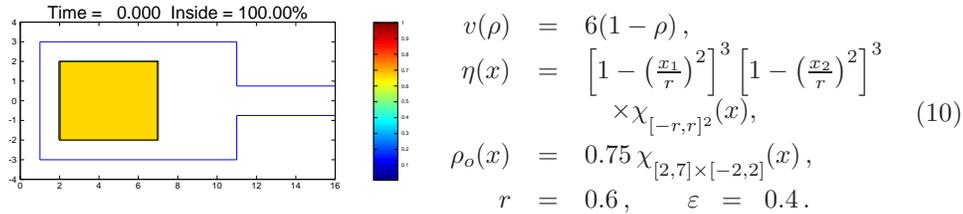}
  \end{minipage}%
  \begin{minipage}[c]{0.57\linewidth}
    \begin{equation}
      \label{eq:Evacuation}
      \begin{array}{@{}rcl@{}}
        v (\rho) & = & 6 (1-\rho)\,,
        \\
        \eta (x)
        & = &
        \left[1-\left(\frac{x_1}{r}\right)^2\right]^3
        \left[1-\left(\frac{x_2}{r}\right)^2\right]^3
        \\
        & & \quad \times
        \caratt{[-r,r]^2} (x),
        \\
        \rho_o (x) & = & 0.75 \, \caratt{[2,7] \times [-2,2]} (x)\,,
        \\
        r & = & 0.6\,,
        \qquad
        \epsilon \;\; = \;\; 0.4\,.
      \end{array}
    \end{equation}
  \end{minipage}
  \caption{Initial datum and room geometry, without
    obstacles.\label{fig:Initial}}
\end{figure}
We consider a room with an exit, as in Figure~\ref{fig:Initial}. The
vector $\vec{\nu} = \vec{\nu} (x)$ is the unit vector tangent at $x$
to the geodesic connecting $x$ to the exit and $\mathcal{I} (\rho) = -
\epsilon \, \left. \left(\nabla (\rho*\eta)\right) \middle/
  \sqrt{1+\norma{\nabla(\rho*\eta)}^2} \right.$,
see~(\ref{eq:Evacuation}).

\begin{figure}[htpb]
  \centering
  \includegraphics[width=0.32\textwidth, trim=75 110 10
  120]{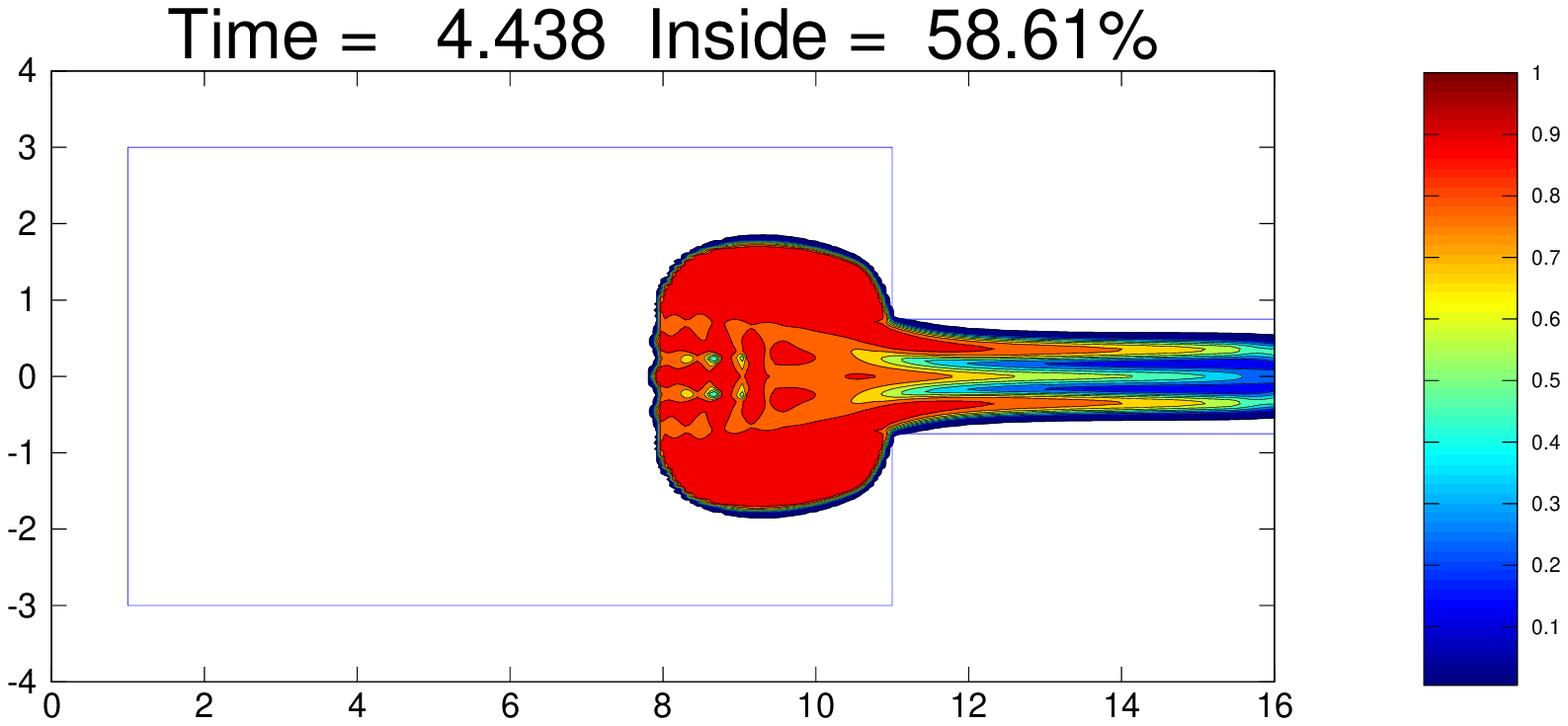}%
  \includegraphics[width=0.32\textwidth, trim=75 110 10
  120]{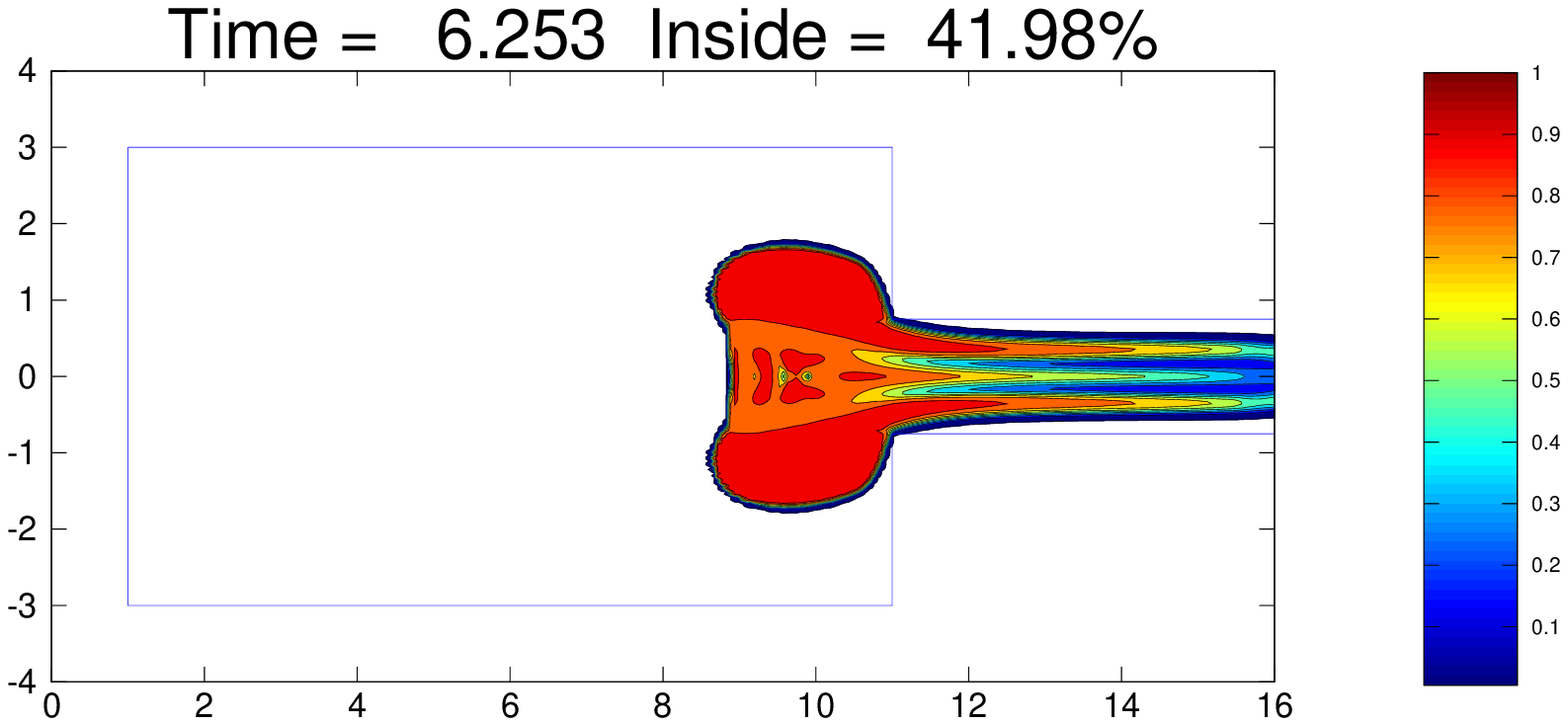}%
  \includegraphics[width=0.32\textwidth, trim=75 110 10 120]{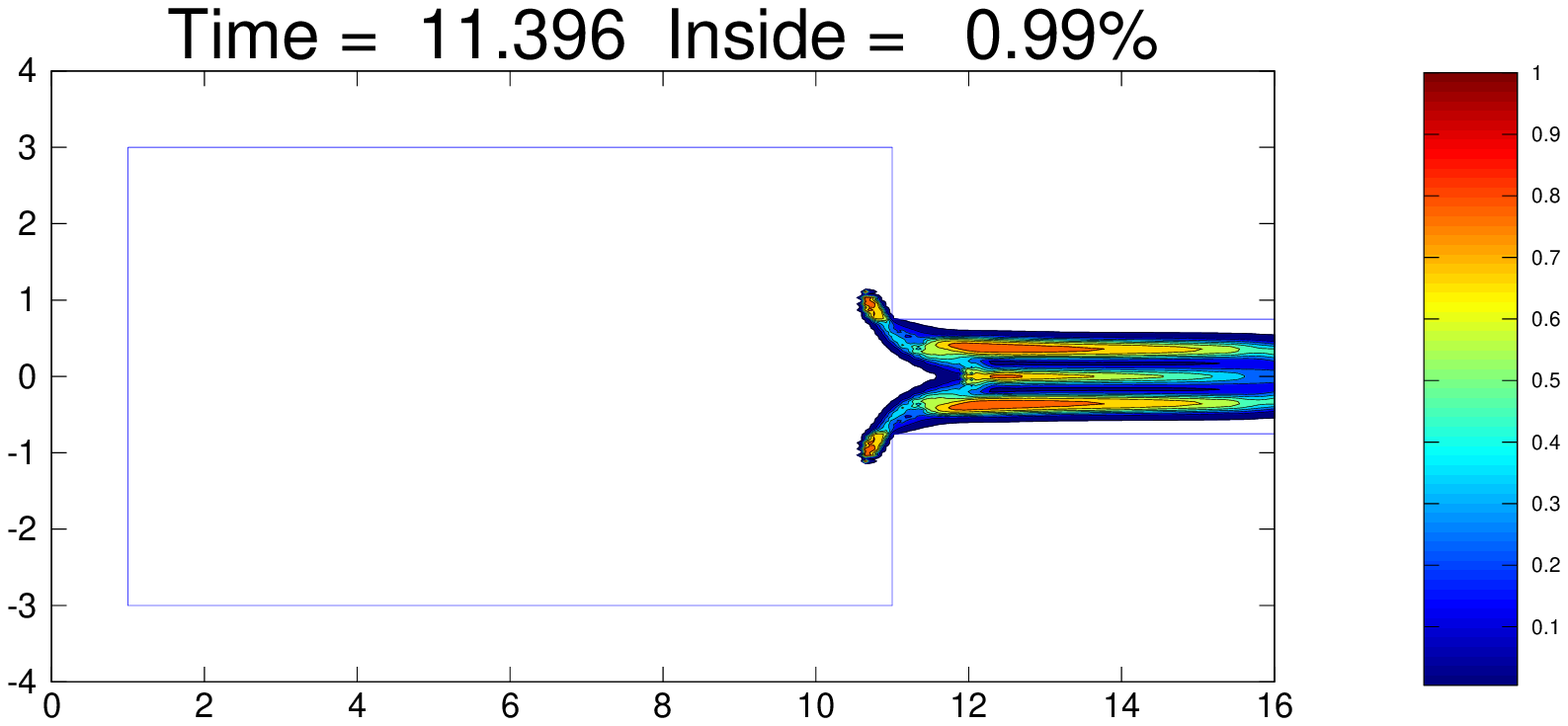}\\
  \includegraphics[width=0.32\textwidth, trim=75 120 10
  110]{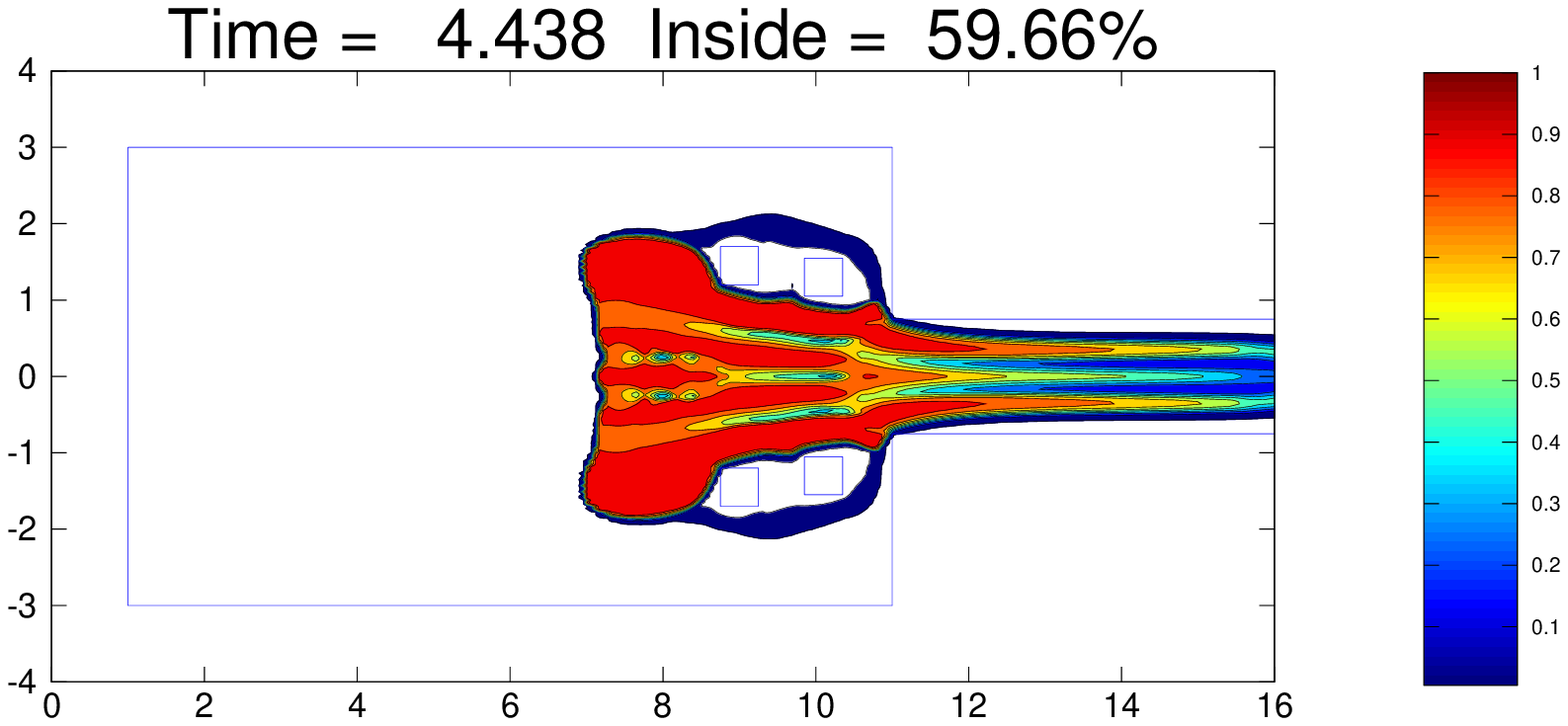}%
  \includegraphics[width=0.32\textwidth, trim=75 120 10
  110]{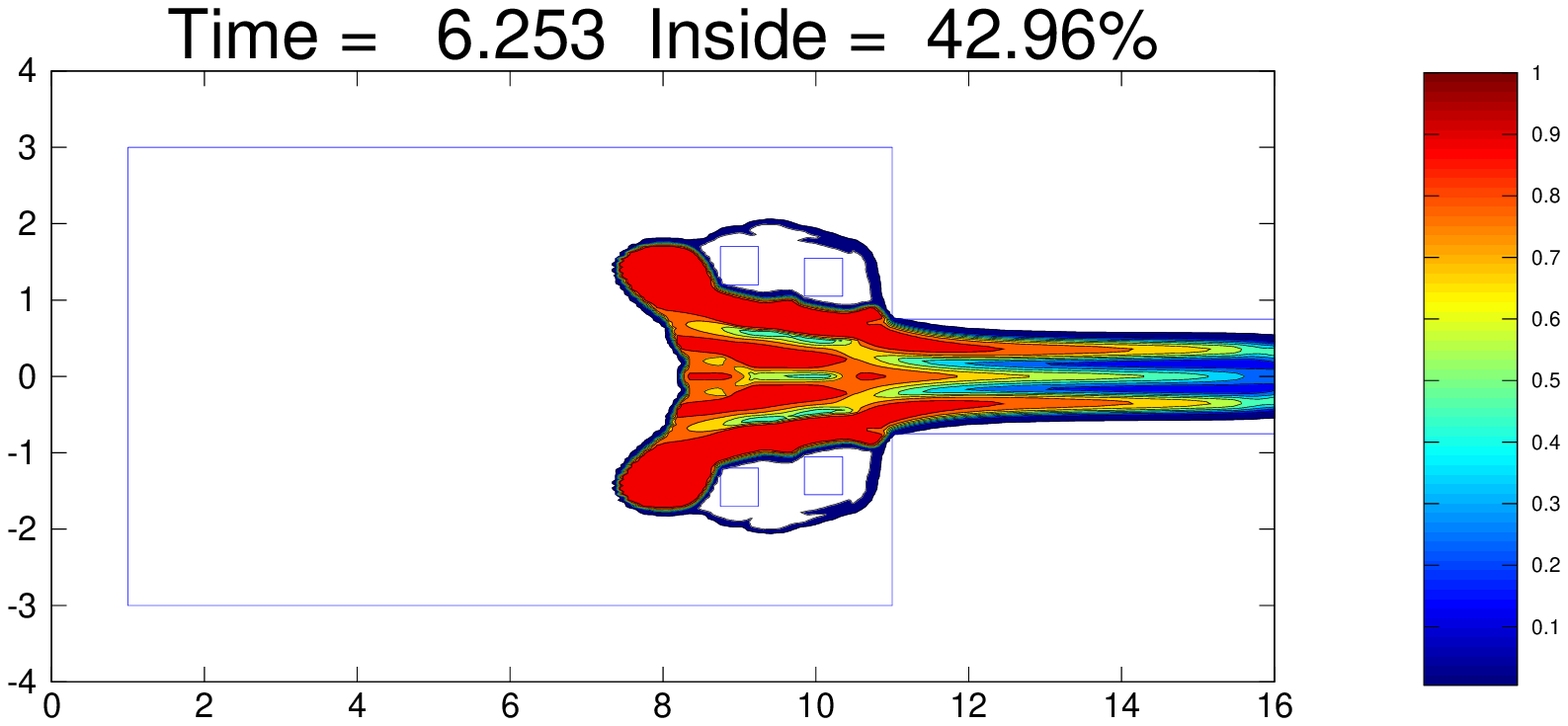}%
  \includegraphics[width=0.32\textwidth, trim=75 120 10 110]{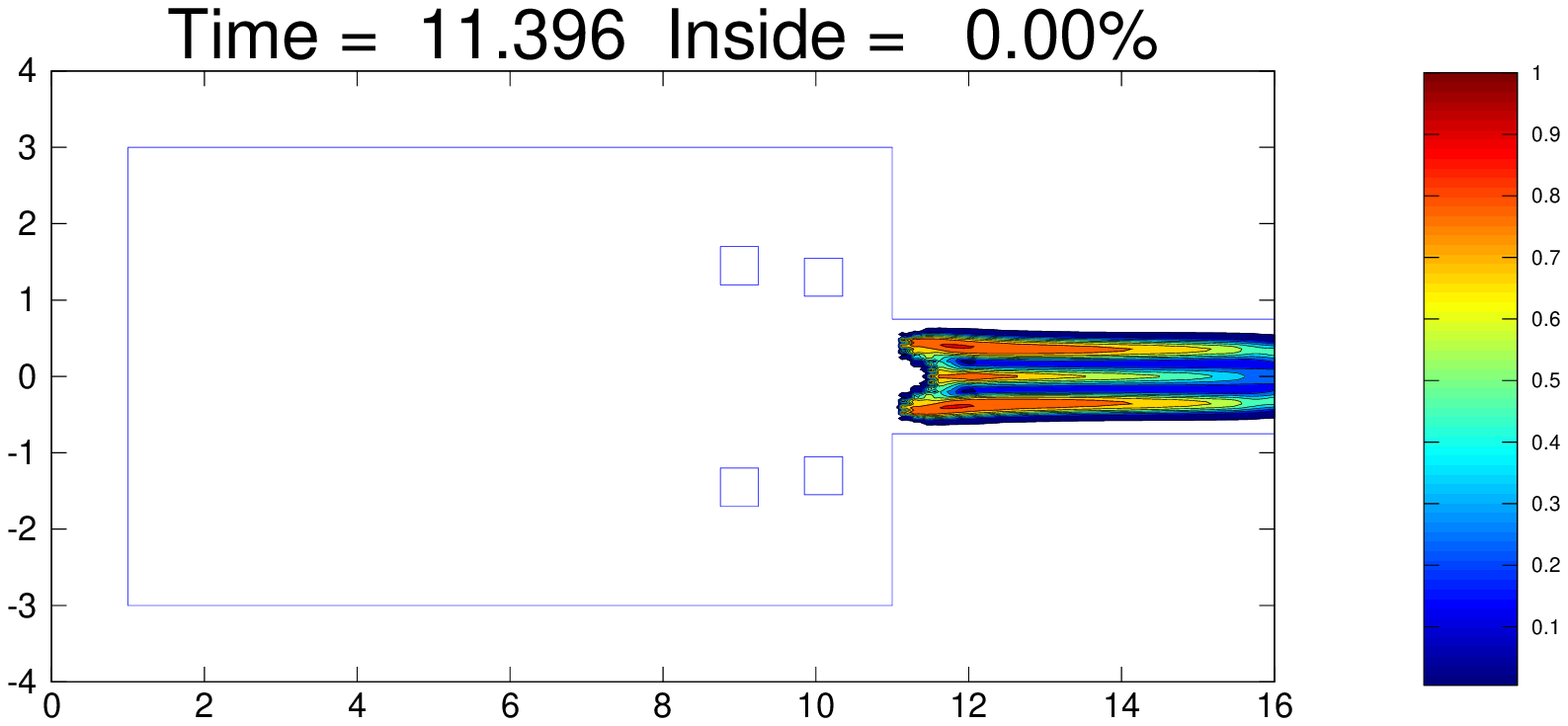}\\
  \caption{Solution
    to~\eqref{eq:1}--\eqref{eq:6}--\eqref{eq:Evacuation} with
    $\epsilon=0.2$, at times $t=4.438$, $6.253$, $11.396$. On the
    first line, no obstacle is present. On the second line, $4$
    columns direct the crowd flow. The exit time in the latter case is
    \emph{shorter} than in the former one,
    see~\cite{ColomboGaravelloLecureux}.\label{fig:Braess}}
\end{figure}
The careful positioning of obstacles as in the second line of
Figure~\ref{fig:Braess} diminishes the size of the congested region
and, with the chosen initial datum, gives an exit time lower than that
with the room free from any obstacle, see Figure~\ref{fig:Braess}.

\section{Individuals Driving a Population}
\label{sec:D}

% Here, briefly quote~\cite{ColomboMercier},
% then~\cite{ColomboPogodaev1, ColomboPogodaev2}, figures.

We finally introduce a model describing the situation in which a
discrete set of isolated individuals interacts with a continuum crowd.
Examples can be a (group of) predator(s) running after their preys,
shepherd dogs driving a herd of sheep, or a leader attracting a group
of followers to a given region.  Let $\rho \in \rpic$ be the
population density and $p \equiv (p_1, \ldots, p_k) \in \reali^{2k}$
be the positions of the $k$
individuals. Following~\cite{ColomboMercier}, the interaction is
described by
\begin{equation}
  \label{eq:groupindiv}
  \left\{
    \begin{array}{l}
      \partial_t \rho
      +
      \div  \left(\rho \, V \! \left(t, x, \rho (t,x), p(t)\right)\right)
      =
      0\,,
      \\
      \dot p
      =
      \phi\left( t, p (t), \rho(t) \right)\,.
    \end{array}
  \right.
  % \qquad
  % \begin{array}{rcl}
  %   (t, x )& \in & \reali^+ \times \reali^{N}\,.
  % \end{array}
\end{equation}
Here, $\phi$ is typically nonlocal, meaning that the individuals $p$
react to averages of quantities depending on $\rho$. The well
posedness of~\eqref{eq:groupindiv} is proved
in~\cite[Theorem~2.2]{ColomboMercier}, by means of Kru\v zkov theory,
the estimates in~\cite{ColomboMercierRosini} and tools from the
stability of ordinary differential equations.

As a first illustrating example, assume that the vector $p\in
\reali^2$ is the position of a leader (e.g.~a magic piper) and $\rho$
is the density of the followers (e.g.~rats). We are thus lead to
consider~\eqref{eq:groupindiv} with
\begin{equation}
  \label{eq:7}
  \begin{array}{rcl}
    V (t,x,\rho,p)
    & = &
    v(\rho) \, (p - x) \, e^{-\norma{p-x}}
    \\
    \phi (t, \pi, \rho)
    & = &
    \left(1 + (\rho*\eta) \left(p(t) \right) \right) \, \vec{{\psi}}(t)\,.
  \end{array}
\end{equation}
The function $v$ essentially describes the speed of the followers and
is, as usual, a smooth decreasing function vanishing at, say,
$\rho=1$. The follower located at $x$ moves along $p(t)-x$ toward the
leader, with a speed exponentially decreasing with the distance
$\norma{p-x}$ between leader and follower. The speed of the leader
increases with the averaged density $\rho*\eta$, computed at the
leader's position. Indeed, we expect the leader to wait for the
followers to join him when the followers' density around him is
small. The direction $\vec{\psi}$ of the leader is chosen \emph{a
  priori}.
\begin{figure}[!htpb]
  \centering
  \includegraphics[width=0.35\textwidth, trim=75 110 50
  120]{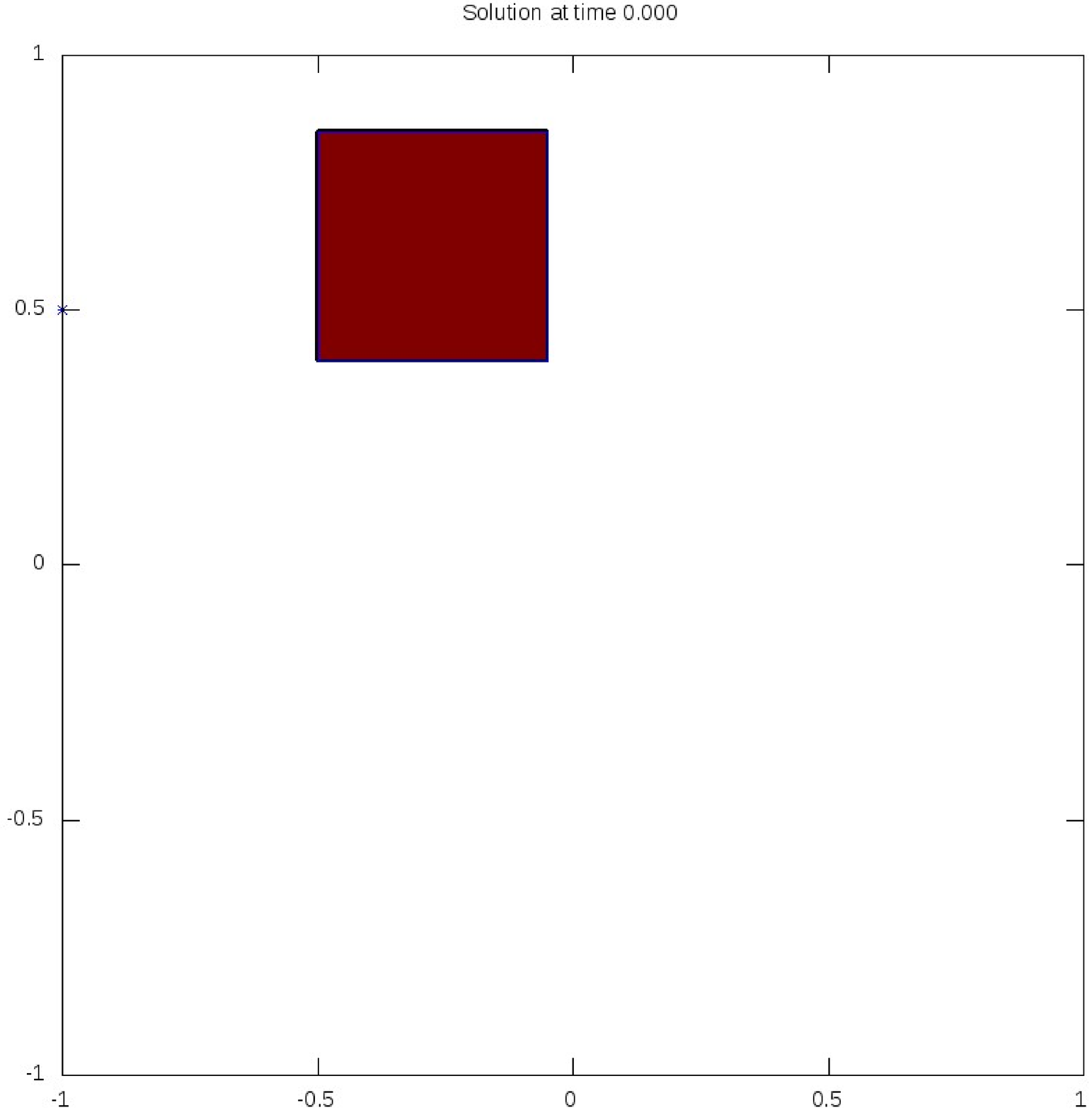}%
  \includegraphics[width=0.35\textwidth, trim=75 110 50
  120]{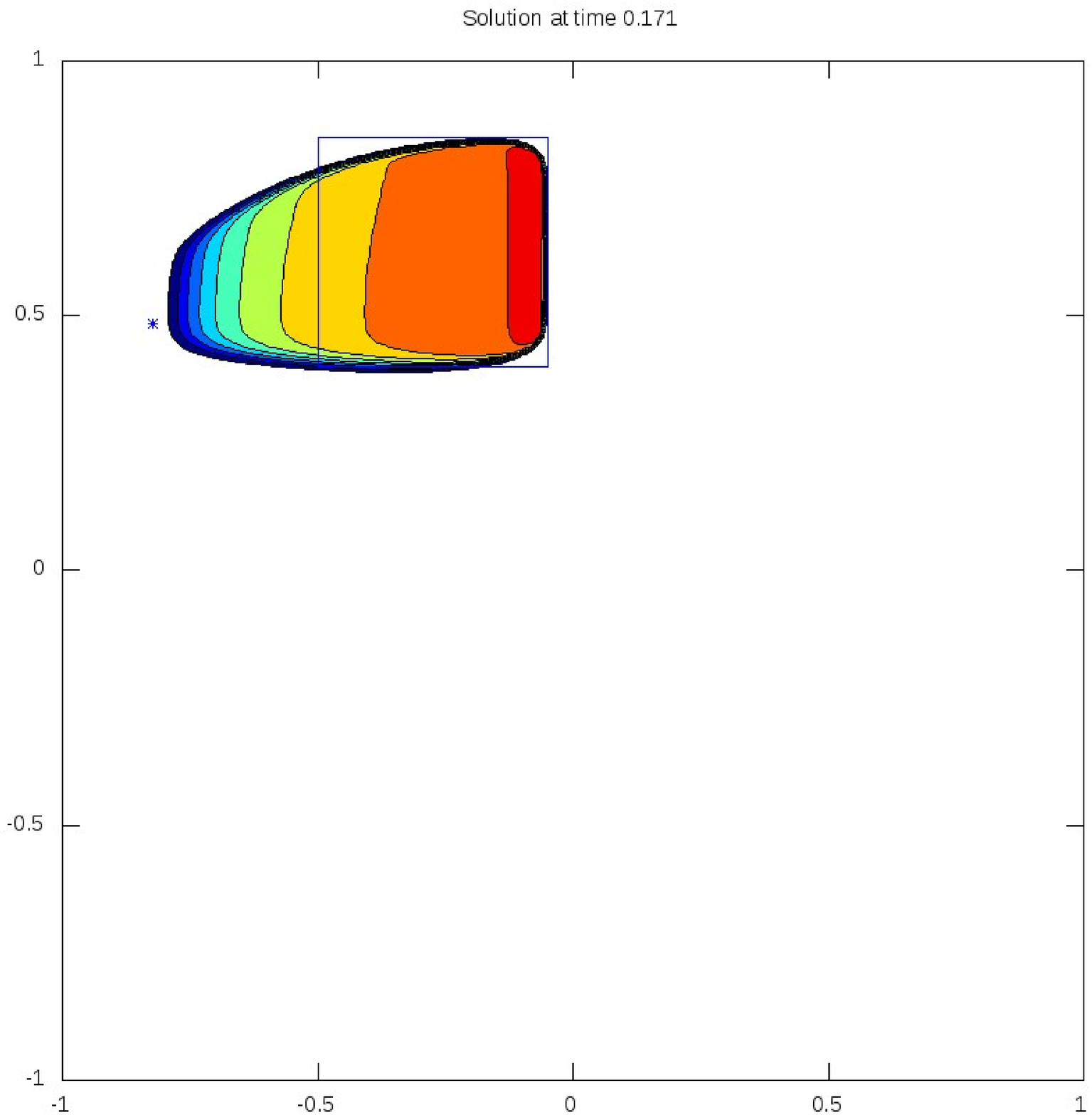}%
  \includegraphics[width=0.35\textwidth, trim=75 110 50
  120]{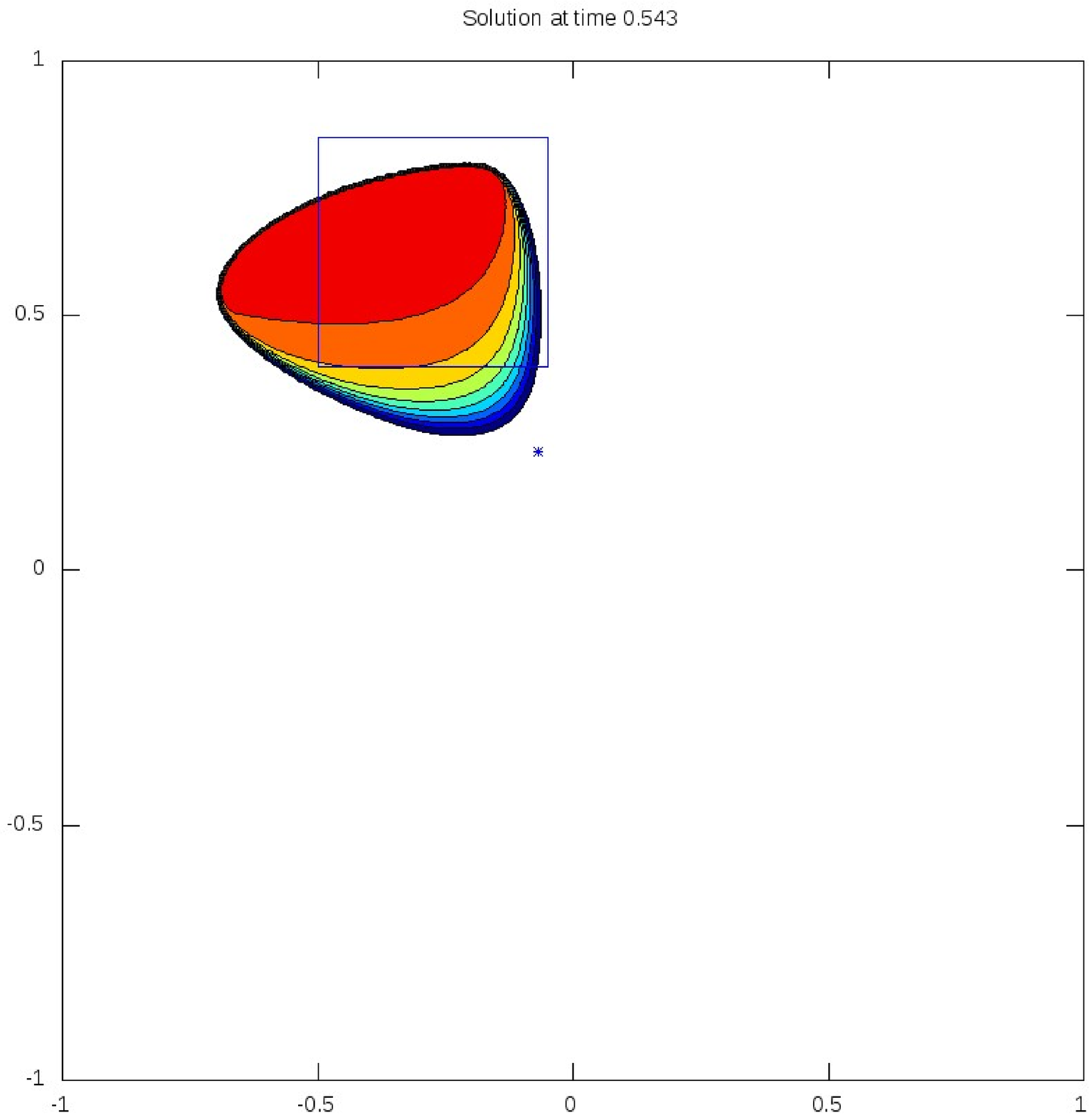}\\
  \includegraphics[width=0.35\textwidth, trim=75 110 50
  120]{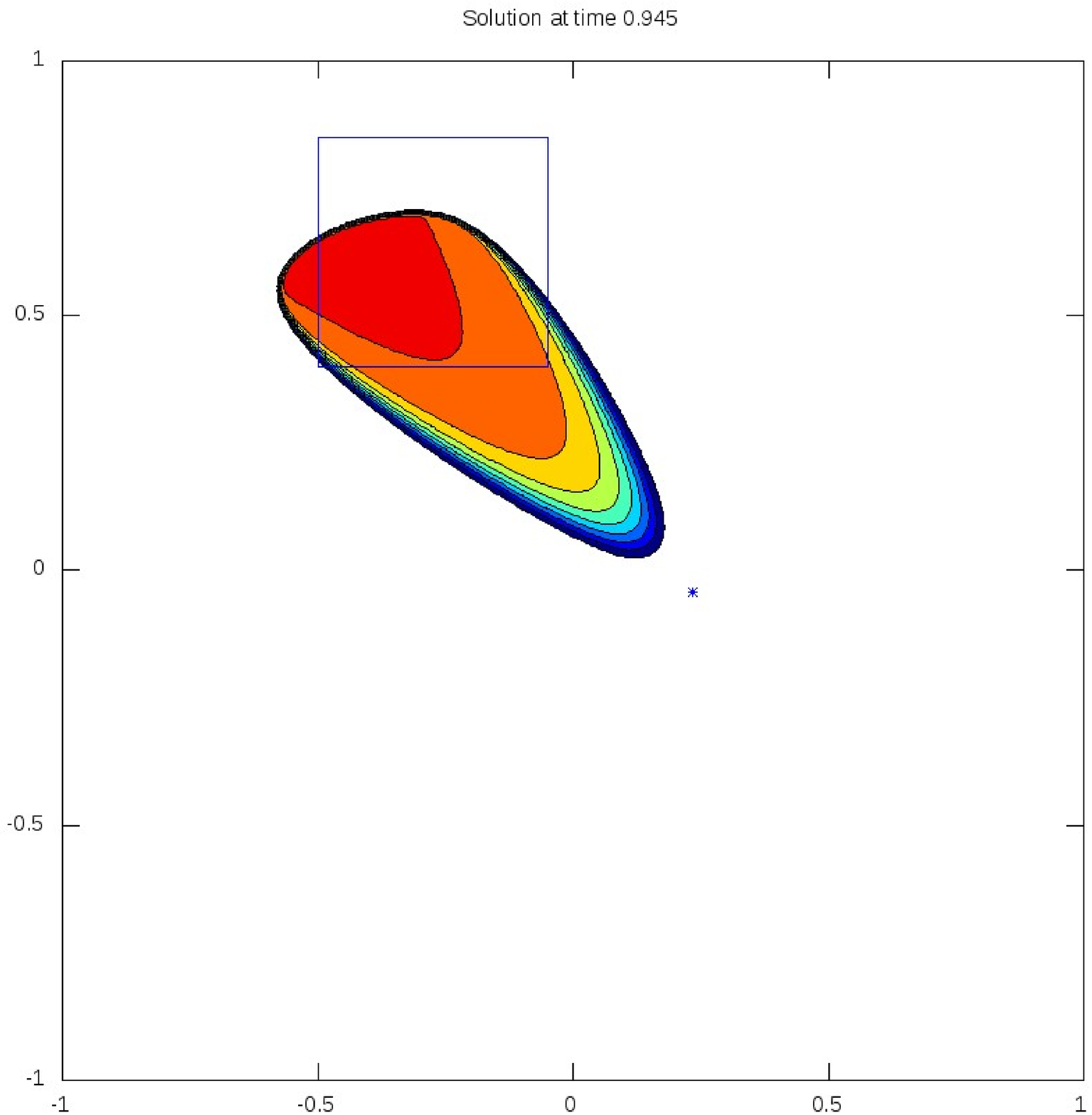}%
  \includegraphics[width=0.35\textwidth, trim=75 110 50
  120]{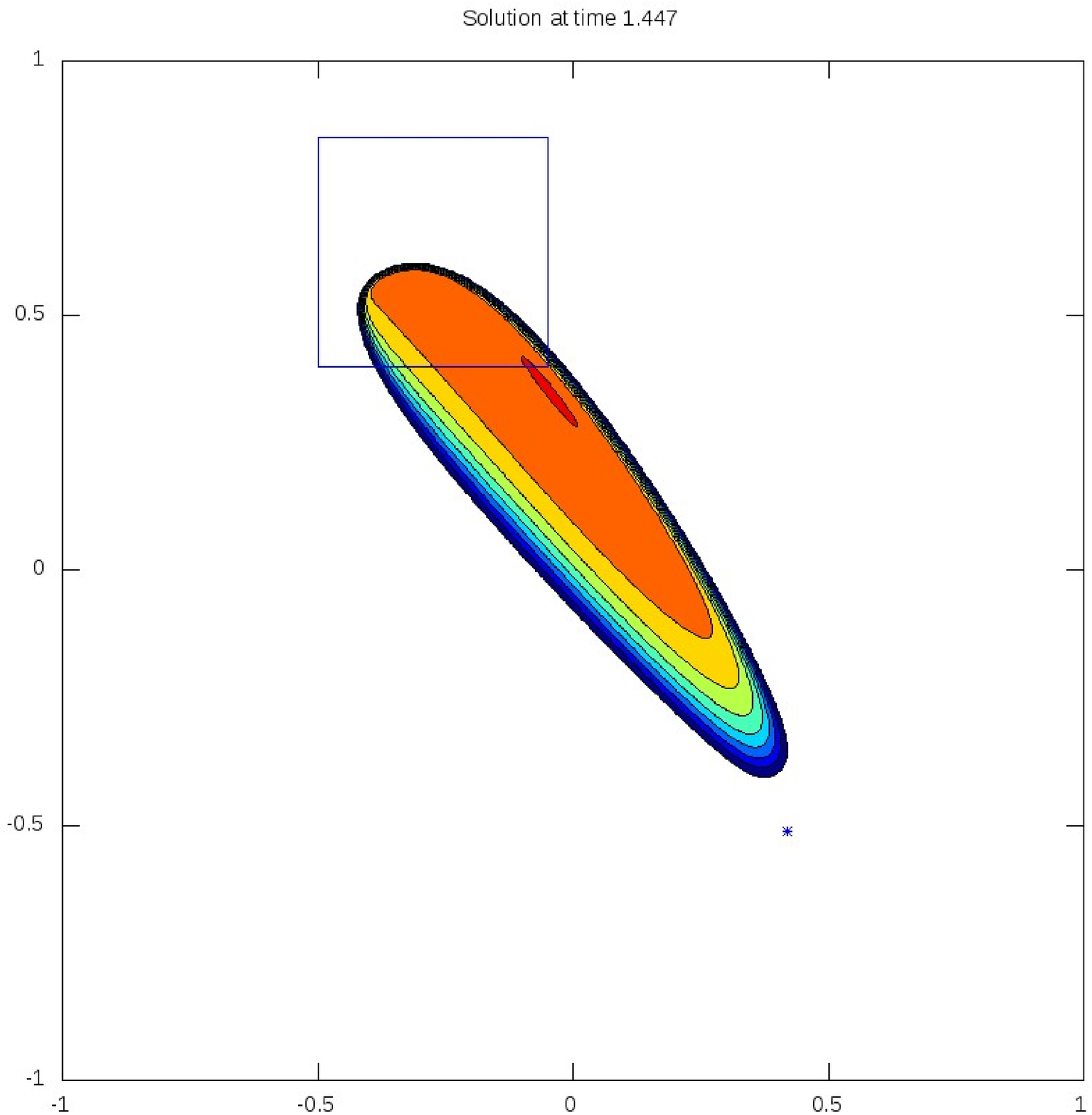}%
  \includegraphics[width=0.35\textwidth, trim=75 110 50
  120]{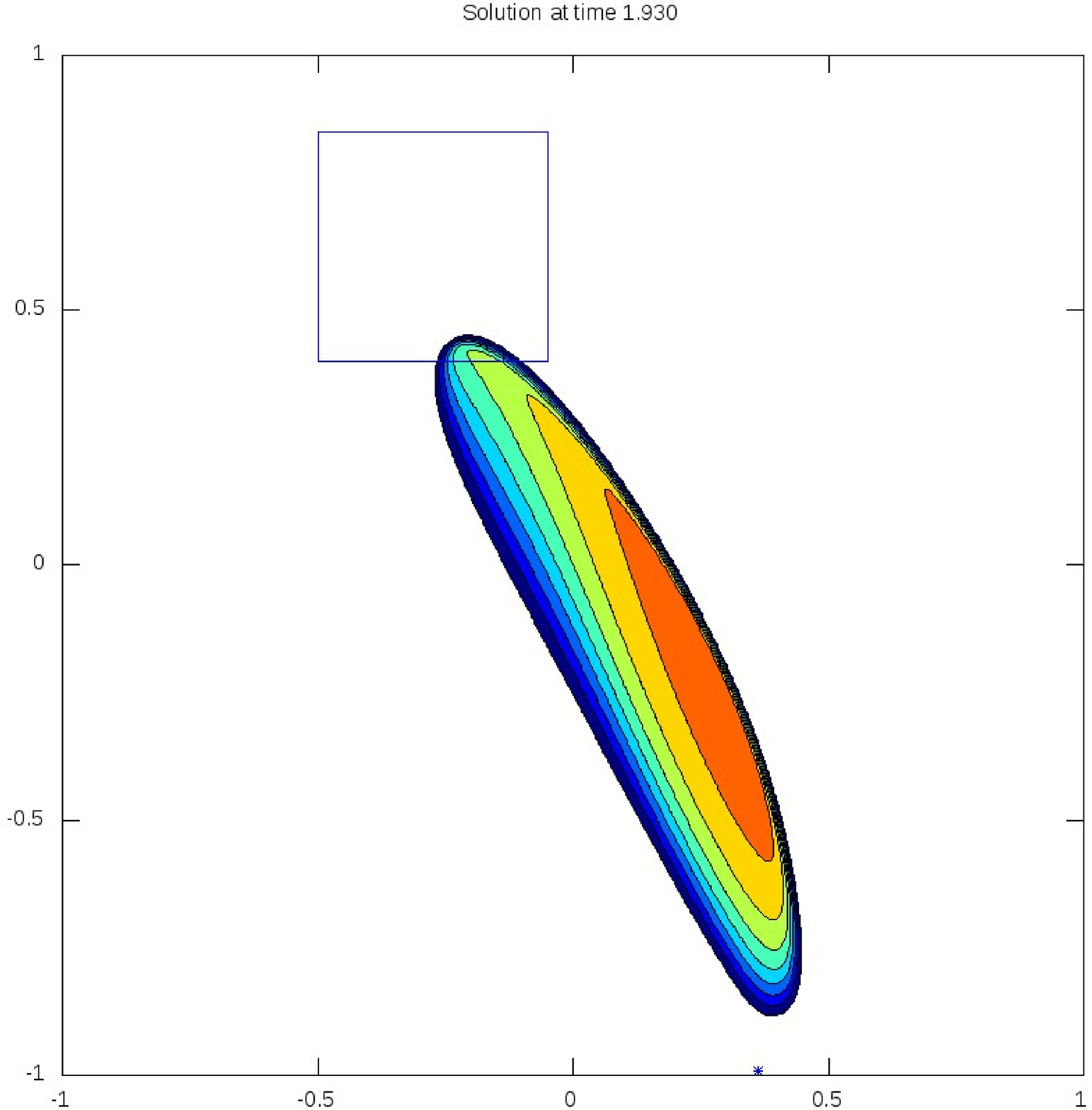}%
  \label{fig:lanes}
  \caption{Solution of~\eqref{eq:groupindiv}--\eqref{eq:7},
    from~\cite{ColomboMercier}.\label{fig:piper}}
\end{figure}
See Figure~\ref{fig:piper} for a numerical integration
of~\eqref{eq:groupindiv}--\eqref{eq:7} and~\cite{ColomboMercier} for
further details.

As a further example, consider $n$ shepherd dogs, located in
$p_i(t)\in \reali^2$ for $i \in \{1, \ldots, n\}$ and a group of sheep
of density $\rho$. The dogs have to confine the sheep within a given
area. We are thus lead to consider~\eqref{eq:groupindiv} with
\begin{equation}
  \label{eq:8}
  \begin{array}{rcl}
    V (t,x,\rho,p)
    & = &
    v(\rho) \, \vec{\nu}(x)
    + \sum_{i=1}^n  (x-p_i)e^{-\norma{p_i-x}}
    \\
    \phi (t, \pi, \rho)
    & = &
    \frac{(\rho *\nabla \eta )^\perp(p_i(t))}{\sqrt{1+\norma{\rho * \nabla \eta(p_i(t)) }^2}}
    \,, \qquad \textrm{for } i\in \{1, \ldots n\} \,.
  \end{array}
\end{equation}
As above, the speed of the sheep is given by the decreasing function
$v(\rho)$ that vanishes in $\rho=1$. The direction of a sheep located
at $x$ is a sum of two terms. The first one is the sheep's preferred
direction $\vec{\nu}(x)$; the second one is the vector $\sum_{i=1}^n
(x-p_i(t)) e^{-\norma{p_i-x}}$ representing the repulsive effects of
the dogs on the sheep. Each dog runs around the flock along the
direction perpendicular to the gradient of the sheep average density.

\section{A Different Approach}
\label{sec:DI}

Following~\cite{BressanZhang, ColomboPogodaev1, ColomboPogodaev2}, we
present another framework to describe the population--individuals
interactions. Initially, the population occupies the compact set $K_o
\subset \reali^2$.  If there are no individuals, the member at $x$ of
the population is free to wander in $\reali^2$, according to the
differential inclusion
\begin{equation}
  \label{eq:di0}
  \dot x \in B (0,c)\,, \qquad x (0) \in K_o \,,
\end{equation}
$c$ being the maximal wandering speed and $B (0,c)$ the closed ball in
$\reali^2$ centered at $0$ with radius $c$. Hence, the population
fills the reachable set of~\eqref{eq:di0}. Introduce now $n$
individuals sited at $\xi \equiv (\xi_1, \xi_2 , \ldots, \xi_n) \in
\reali^{2n}$. Then, the interaction between the individuals and each
population member leads to the modified differential inclusion
\begin{equation}
  \label{eq:di}
  \dot x \in v \left(x,\xi (t)\right) +  B(0,c),\quad x(0)\in K_o,
\end{equation}
where the vector field $v \in \C{0,1} (\reali^2 \times \reali^{2n};
\reali^2)$ is the drift speed due to the attractive or repulsive
effect that each agent has on each member of the population.  Thus,
given the individuals' trajectory $\xi \in
\Cloc{0,1}(\reali^+;\reali^{2n})$, the reachable set
$\mathcal{R}_\xi(K_o,t)$ of~\eqref{eq:di} at time $t$ is the set
occupied by the population at time $t$ under the effect of the agents.
With the present assumptions, $\mathcal{R}_\xi(K_o,t)$ is non-empty
and compact.

If only one agent is present ($n=1$) and $v$ is spherically symmetric,
i.e.,
\begin{equation}
  \label{eq:v}
  v(x,\xi) = \psi(|x-\xi|)(x-\xi)\quad\mbox{for a suitable }\; \psi\colon\reali\to\reali.%\in \C{0,1}(\reali;\reali).
\end{equation}
the next result exhibits a trajectory $\xi$ confining the population
in a given set $K$.

\begin{theorem}{~\cite[Theorem 2.8]{ColomboPogodaev2}}
  \label{thm:pos}
  Let $c > 0$. Fix a bounded $\psi \in \Cloc{1,1} (\reali^n; \reali)$
  and define $v$ as in~\eqref{eq:v}. Assume that there exist positive
  $R_*^-$, $R_*^+$ and $R$ such that
  \begin{equation*}
    \frac{1}{\pi}
    \int_0^\pi
    \psi \left(\sqrt{R^2 + R_*^2 - 2 R_* \, R\, \cos\theta} \right)
    (R_* - R \, \cos\theta) \, \d\theta
    <
    -c
  \end{equation*}
  for all $R_* \in [R_*^-,\, R_*^+]$.  Then, there exists a $\xi \in
  \Cloc{0,1} \left(\reali^+; \partial B(0,R)\right)$ such that,
  calling $K_o$ the region initially occupied by the population,
  \begin{equation*}
    \mbox{if } \quad
    K_o \subseteq B (0, R_*^-)
    \quad \mbox{ then } \quad
    \mathcal{R}_\xi(t, K_o) \subseteq B (0, R_*^+)
    \quad \mbox{ for all } \quad t\geq 0 \,.
  \end{equation*}
\end{theorem}
Note that the confining strategy $t\mapsto \xi(t)$ above is
constructed explicitly, see~\cite[Theorem~2.5]{ColomboPogodaev1}. A
negative result is also available.  Before stating it, recall that for
a measurable function $\phi \colon \reali^+ \to \reali$, its
\emph{non--decreasing rearrangement} is the function $\phi_* \colon
\reali^+ \to \reali$, which is non--decreasing and satisfies
$\mathcal{L}^1 \left( \phi_*^{-1}(\left]-\infty, a\right]) \right) =
\mathcal{L}^1 \left( \phi^{-1}(\left]-\infty, a\right]) \right)$ for
all $a \in \reali$.

\begin{theorem}{~\cite[Theorem 2.7]{ColomboPogodaev2}}
  \label{thm:neg}
  Let $c > 0$. Fix a bounded $\psi \in \Cloc{1,1} (\reali; \reali)$
  and define $v$ as in~\eqref{eq:v}. Let $\phi_*$ be the
  non--decreasing rearrangement of the function
  \begin{equation*}
    \phi(s) =
    \psi' \! \left( \sqrt[2]{\frac{s}{\pi}} \right)
    \sqrt[2]{\frac{s}{\pi}}
    +
    2\, \psi \! \left(\sqrt[2]{\frac{s}{\pi}}\right) \,.
  \end{equation*}
  If the initial set $K_o$ is such that
  \begin{equation*}
    \label{eq:g}
    2 \, c \, \sqrt{\pi\sigma}
    +
    \int_0^\sigma \phi_*(s) \, \d s
    > 0
    \quad \mbox{ for all } \quad
    \sigma \geq \mathcal{L}^2 (K_o)
  \end{equation*}
  then, for every $\xi \in \Cloc{0,1}(\reali^+;\reali^{kn})$, the
  measure $\mathcal{L}^2\left(\mathcal{R}_\xi (t, K_o)\right)$ of the
  reachable set $\mathcal{R}_\xi (t, K_o)$ of~\eqref{eq:di} increases
  unboundedly in time, so that no confinement is possible.
\end{theorem}
We refer to the cited references for the statement of these results in
arbitrary space dimension. Theorem~\ref{thm:neg} holds also in the
case of several individuals, each acting as in~\eqref{eq:v}
(see~\cite{ColomboPogodaev2}).

\bibliography{ColomboGaravelloMercierPogodaev}

\begin{thebibliography}{10}

\bibitem{AmadoriDiFrancesco}
D.~Amadori and M.~Di~Francesco.
\newblock The one-dimensional {H}ughes model for pedestrian flow:
  {R}iemann--type solutions.
\newblock {\em Acta Mathematica Scientia}, 32(1):259--196, 2011.

\bibitem{BellomoDogbe2011}
N.~Bellomo and C.~Dogbe.
\newblock On the modeling of traffic and crowds: a survey of models,
  speculations, and perspectives.
\newblock {\em SIAM Rev.}, 53(3):409--463, 2011.

\bibitem{BraessParadox}
D.~Braess.
\newblock \"{U}ber ein {P}aradoxon aus der {V}erkehrsplanung.
\newblock {\em Unternehmensforschung}, 12:258--268, 1968.

\bibitem{BressanLectureNotes}
A.~Bressan.
\newblock {\em Hyperbolic systems of conservation laws}, volume~20 of {\em
  Oxford Lecture Series in Mathematics and its Applications}.
\newblock Oxford University Press, Oxford, 2000.
\newblock The one-dimensional Cauchy problem.

\bibitem{BressanFire}
A.~Bressan.
\newblock Differential inclusions and the control of forest fires.
\newblock {\em J. Differential Equations}, 243(2):179--207, 2007.

\bibitem{BressanZhang}
A.~Bressan and D.~Zhang.
\newblock Control problems for a class of set valued evolutions.
\newblock {\em Set-Valued and Variational Analysis}, pages 1--21, 2012.
\newblock 10.1007/s11228-012-0204-5.

\bibitem{ColomboGaravelloLecureux}
R.~M. Colombo, M.~Garavello, and M.~L{\'e}cureux-Mercier.
\newblock Non-local crowd dynamics.
\newblock {\em C. R. Math. Acad. Sci. Paris}, 349(13-14):769--772, 2011.

\bibitem{ColomboGaravelloLecureuxM3AS}
R.~M. Colombo, M.~Garavello, and M.~L{\'e}cureux-Mercier.
\newblock A class of non-local models for pedestrian traffic.
\newblock {\em Mathematical Models and Methods in the Applied Sciences}, 22(4),
  2012.

\bibitem{ColomboHertyMercier}
R.~M. Colombo, M.~Herty, and M.~Mercier.
\newblock Control of the continuity equation with a non local flow.
\newblock {\em ESAIM Control Optim. Calc. Var.}, 17(2):353--379, 2011.

\bibitem{ColomboLecureuxPerDafermos}
R.~M. Colombo and M.~L{\'e}cureux-Mercier.
\newblock Nonlocal crowd dynamics models for several populations.
\newblock {\em Acta Mathematica Scientia}, 32(1):177--196, 2011.

\bibitem{ColomboMercier}
R.~M. Colombo and M.~Mercier.
\newblock An analytical framework to describe the interactions between
  individuals and a continuum.
\newblock {\em Journal of Nonlinear Science}, 22(1):39--61, 2012.

\bibitem{ColomboMercierRosini}
R.~M. Colombo, M.~Mercier, and M.~D. Rosini.
\newblock Stability and total variation estimates on general scalar balance
  laws.
\newblock {\em Commun. Math. Sci.}, 7(1):37--65, 2009.

\bibitem{ColomboPogodaev1}
R.~M. Colombo and N.~Pogodaev.
\newblock Confinement strategies in a model for the interaction between
  individuals and a continuum.
\newblock {\em SIAM J. Appl. Dyn. Syst.}, 11(2):741--770, 2012.

\bibitem{ColomboPogodaev2}
R.~M. Colombo and N.~Pogodaev.
\newblock On the control of moving sets: Positive and negative confinement
  results.
\newblock {\em Preprint}, 2012.

\bibitem{ColomboRosini2005}
R.~M. Colombo and M.~D. Rosini.
\newblock Pedestrian flows and non-classical shocks.
\newblock {\em Math. Methods Appl. Sci.}, 28(13):1553--1567, 2005.

\bibitem{ColomboRosini2009}
R.~M. Colombo and M.~D. Rosini.
\newblock Existence of nonclassical solutions in a pedestrian flow model.
\newblock {\em Nonlinear Anal. Real World Appl.}, 10(5):2716--2728, 2009.

\bibitem{CrippaMercier}
G.~Crippa and M.~L\'ecureux-Mercier.
\newblock Existence and uniqueness of measure solutions for a system of
  continuity equations with non-local flow.
\newblock {\em NODEA}, 2012.

\bibitem{CristianiPiccoliTosin}
E.~Cristiani, B.~Piccoli, and A.~Tosin.
\newblock Multiscale modeling of granular flows with application to crowd
  dynamics.
\newblock {\em Multiscale Model. Simul.}, 9(1):155--182, 2011.

\bibitem{DafermosBook}
C.~M. Dafermos.
\newblock {\em Hyperbolic conservation laws in continuum physics}, volume 325
  of {\em Grundlehren der Mathematischen Wissenschaften [Fundamental Principles
  of Mathematical Sciences]}.
\newblock Springer-Verlag, Berlin, third edition, 2010.

\bibitem{DiFrancescoMarkowichWolfram}
M.~Di~Francesco, P.~A. Markowich, J.-F. Pietschmann, and M.-T. Wolfram.
\newblock On the {H}ughes' model for pedestrian flow: the one-dimensional case.
\newblock {\em J. Differential Equations}, 250(3):1334--1362, 2011.

\bibitem{GoatinRosini}
N.~El-Khatib, P.~Goatin, and M.~D. Rosini.
\newblock On entropy weak solutions of {H}ughes’ model for pedestrian motion.
\newblock {\em ZAMP}, 2012.
\newblock To appear.

\bibitem{HelbingJohanssonZein}
D.~Helbing, A.~Johansson, and H.~Z. Al-Abideen.
\newblock Dynamics of crowd disasters: An empirical study.
\newblock {\em Phys. Rev. E}, 75:046109, Apr 2007.

\bibitem{HelbingEtAlii2001}
D.~Helbing, P.~Moln\'ar, I.~Farkas, and K.~Bolay.
\newblock Self-organizing pedestrian movement.
\newblock {\em Environment and Planning B: Planning and Design}, 28:361--383,
  2001.

\bibitem{Kruzkov}
S.~N. Kru{\v{z}}hkov.
\newblock First order quasilinear equations with several independent variables.
\newblock {\em Mat. Sb. (N.S.)}, 81 (123):228--255, 1970.

\bibitem{MauryEtAl}
B.~Maury, A.~Roudneff-Chupin, and F.~Santambrogio.
\newblock A macroscopic crowd motion model of gradient flow type.
\newblock {\em Math. Models Methods Appl. Sci.}, 20(10):1787--1821, 2010.

\bibitem{PiccoliTosin}
B.~Piccoli and A.~Tosin.
\newblock Time-evolving measures and macroscopic modeling of pedestrian flow.
\newblock {\em Arch. Ration. Mech. Anal.}, 199(3):707--738, 2011.

\end{thebibliography}

\bibliographystyle{abbrv}

\end{document}